\documentclass[11pt]{article}

\usepackage{epsfig,amsmath,latexsym}
\usepackage{amsfonts}

\newtheorem{theorem}{Theorem}[section]

\newtheorem{defi}{Definition}[section]
\newtheorem{lemma}[theorem]{Lemma}

\def\slfrac#1#2{\hbox{\kern.1em %
 \raise.5ex\hbox{\the\scriptfont0 #1}\kern-.11em %
 /\kern-.15em\lower.25ex\hbox{\the\scriptfont0 #2}}}

\newcommand{\pf}{\noindent{\bf Proof.~}}
\newcommand{\beq}{\begin{eqnarray}}
\newcommand{\eeq}{\end{eqnarray}}
\newcommand{\beql}[1]{\begin{eqnarray}\label{#1}}
\newcommand{\beqs}{\begin{eqnarray*}}
\newcommand{\eeqs}{\end{eqnarray*}}
\newcommand{\eqn}[1]{(\ref{#1})}

\newcommand{\rr}{{\mathbb R}}
\newcommand{\sss}{{\mathbb S}}
\newcommand{\zz}{{\mathbb Z}}

\newcommand{\hh}{{\mathbb H}}

\newcommand{\bb}{{\mathbf b}}
\newcommand{\bc}{{\mathbf c}}

\newcommand{\bff}{{\mathbf f}}
\newcommand{\bg}{{\mathbf g}}
\newcommand{\bh}{{\mathbf h}}
\newcommand{\bm}{{\mathbf m}}

\newcommand{\bu}{{\mathbf u}}

\newcommand{\bw}{{\mathbf w}}
\newcommand{\bx}{{\mathbf x}}
\newcommand{\by}{{\mathbf y}}
  
\newcommand{\bo}{{\mathbf 1}}

\newcommand{\bA}{{\mathbf A}}
\newcommand{\bB}{{\mathbf B}}

\newcommand{\bF}{{\mathbf F}}
\newcommand{\bG}{{\mathbf G}}

\newcommand{\bI}{{\mathbf I}}
\newcommand{\bJ}{{\mathbf J}}

\newcommand{\bM}{{\mathbf M}}
\newcommand{\bN}{{\mathbf N}}

\newcommand{\bQ}{{\mathbf Q}}

\newcommand{\bW}{{\mathbf W}}

\makeatletter
\def\@sect#1#2#3#4#5#6[#7]#8{\ifnum #2>\c@secnumdepth
     \def\@svsec{}\else
     \refstepcounter{#1}\edef\@svsec{\csname the#1\endcsname.\hskip .75em }\fi
     \@tempskipa #5\relax
      \ifdim \@tempskipa>\z@
        \begingroup #6\relax
          \@hangfrom{\hskip #3\relax\@svsec}{\interlinepenalty \@M #8\par}%
        \endgroup
       \csname #1mark\endcsname{#7}\addcontentsline
         {toc}{#1}{\ifnum #2>\c@secnumdepth \else
                      \protect\numberline{\csname the#1\endcsname}\fi
                    #7}\else
        \def\@svsechd{#6\hskip #3\@svsec #8\csname #1mark\endcsname
                      {#7}\addcontentsline
                           {toc}{#1}{\ifnum #2>\c@secnumdepth \else
                             \protect\numberline{\csname the#1\endcsname}\fi
                       #7}}\fi
     \@xsect{#5}}
\def\@begintheorem#1#2{\it \trivlist \item[\hskip \labelsep{\bf #1\ #2.}]}

\def\plain{plain}\ifx\fmtname\plain\csname fi\endcsname
     
     \input docstrip
     \preamble

     Do not distribute the stripped version of this file.
     The checksum in the header refers to the documented version.

     \endpreamble
     \generateFile{here.sty}{t}{\from{here.doc}{}}
     \endinput
\fi
\ifcat a\noexpand @\let\next\relax\else\def\next{%
    \documentstyle[here,doc]{article}\MakePercentIgnore}\fi\next
\ifx\@Hxfloat\@Hundef\else\expandafter\endinput\fi
\let\@Hxfloat\@xfloat
\def\@xfloat#1[{\@ifnextchar{H}{\@HHfloat{#1}[}{\@Hxfloat{#1}[}}
\def\@HHfloat#1[H]{%
\expandafter\let\csname end#1\endcsname\end@Hfloat
\vskip\intextsep\vbox\bgroup\def\@captype{#1}\parindent\z@
\ignorespaces}
\def\end@Hfloat{\egroup\vskip \intextsep}

\makeatother

\catcode`\@=11

\catcode`@=12

\begin{document}

\begin{center}
{\Large {\bf Beyond the Descartes Circle Theorem}}\\
\vspace*{.2\baselineskip}
{\em Jeffrey C. Lagarias}\\
{\em Colin L. Mallows} \\
\vspace*{.2\baselineskip}
{\em Allan R. Wilks} \\
\vspace*{.2\baselineskip}
AT\&T Labs, 
Florham Park, NJ 07932-0971 \\
\vspace*{1.5\baselineskip}
(January 8, 2001)\\
\vspace*{1.5\baselineskip}
{\bf ABSTRACT}
\end{center}
The Descartes circle theorem states that if four circles
are mutually tangent in the plane, with disjoint interiors,
then their curvatures (or ``bends'')
$b_i = \frac {1}{r_i}$ satisfy the 
relation 
$(b_1 +b_2 + b_3 + b_4)^2 = 2 (b_1^2 + b_2^2 + b_3^2 + b_4^2).$
We show that similar relations hold involving the
centers of the four circles in such a configuration, 
coordinatized as
complex numbers, yielding a complex Descartes Theorem. 
These relations have elegant matrix generalizations 
to the $n$-dimensional case, in  each of
Euclidean, spherical, and hyperbolic geometries. These include
analogues of the Descartes circle theorem for spherical and 
hyperbolic space. \\



%
%

\section{Introduction}
\setcounter{equation}{0}

\begin{figure}[htbp]
\centerline{\epsfxsize=2.0in \epsfbox{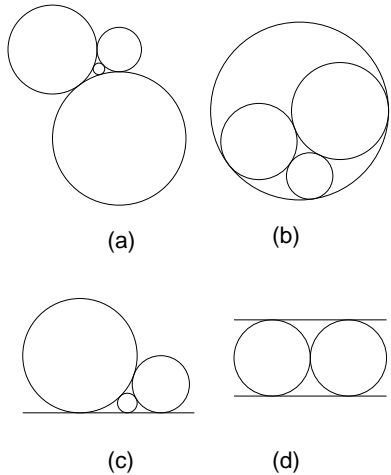}}
\caption{Descartes configurations}\label{fig1}
\end{figure}

We  call a configuration of four mutually tangent circles, 
in which no three circles have a common tangent,
a ``Descartes configuration''. The possible
arrangements~\footnote{An arrangement of concentric circles
having a common tangent is not a Descartes configuration.}
 of such  configurations appear
in Figure~\ref{fig1}, where we
 allow certain ``degenerate'' arrangements
where some of the circles are straight lines.
Suppose the radii of the circles are
$r_1 , r_2 , r_3 , r_4$.  The reciprocals of these are the 
curvatures (or ``bends'')   $b_j = 1 / r_j $.
A straight line is assigned infinite radius; then the ``bend'' is zero.

In 1643 Rene Descartes~\cite[pp. 45--50]{Des01}, in a letter
to Princess Elizabeth of Bohemia, stated a 
relation connecting the four
radii. This relation can be written as a quadratic equation connecting 
the four curvatures, namely:

\begin{theorem}[Descartes Circle Theorem]~\label{th11}
For a Descartes configuration of four mutually tangent circles, their
curvatures satisfy
\beql{101}
\sum_{j = 1}^4 {b_j}^2 = \frac{1}{2}  ( \sum_{j=1}^4 b_j )^2.
\eeq
\end{theorem}

Descartes considered
only the configuration (a) in Figure~\ref{fig1}.
He did not state the result in this form, but 
gave a more complicated relation
algebraically equivalent to \eqn{101},
and his sketched proof is incomplete.
In 1826 Jakob Steiner ~\cite[pp. 61--63]{St26}
independently found the result and gave
a complete proof.  Another
independent rediscovery 
with a complete proof was given in 1842 by H. Beecroft~\cite{Bee42},
and is described in Coxeter~\cite{Cox68}. Many other proofs
have been discovered (and rediscovered), some of which
appear in Pedoe~\cite{Pe67}.

The Descartes circle theorem applies
to all Descartes configurations
of types (a)- (d), provided
we define the curvatures to have appropriate signs,
as follows.
An {\em oriented circle} is a circle together with
an assigned direction of  unit normal vector, 
which can point inward or outward. If it
has radius $r$ then its {\em oriented radius}
is $r$ for an inward pointing normal and $-r$ for
an outward pointing normal.
Its {\em oriented curvature},
(or ``signed curvature'') 
is $\frac{1}{r}$ for an inward pointing normal  and
$- \frac{1}{r}$ for an outward pointing normal.
By convention, the  {\em interior} of an oriented
circle is its interior for an inward pointing normal
and its exterior for an outward pointing normal.
We define an {\em oriented
Descartes configuration} to be a Descartes configuration with the
circles having orientations which are compatible in the following sense:
either (i) the interiors of all four oriented circles are disjoint, or
(ii) the interiors are disjoint when all orientations are reversed.
Each Descartes
configuration has exactly two compatible orientations in this sense, one
obtained from the other by reversing all orientations~\footnote
{The 
{\em inward pointing orientation} of a Descartes configuration
is the one in which at least
two oriented curvatures are strictly positive; the 
{\em outward pointing orientation}
is one in which at least two curvatures are strictly negative.}.
With these definitions,
the Descartes Circle Theorem remains valid  for all oriented Descartes
configurations, using oriented curvatures.

In 1936 Sir Frederick Soddy  (who earned a 1921 Nobel prize for discovering 
isotopes) published in Nature~\cite{Sod36} a poem entitled 
``The Kiss Precise'' in which he reported the result above 
and a generalization to three dimensions. The following year Thorold 
Gossett~\cite{Gos37}
contributed another stanza giving the general $n$-dimensional result.
We extend the definition of a {\em Descartes configuration} to consist 
of $n+2$ mutually
tangent $(n-1)$-spheres in $\rr^n$ in which all pairs of tangent 
$(n-1)$-spheres have
distinct points of tangency, and orientation is done as above.

\begin{theorem}[Soddy-Gossett Theorem]~\label{th12}
Given an oriented  Descartes configuration in $\rr^n$,  
if we let
$b_j = 1 / r_j$ be the oriented  curvatures of the $n+2$ mutually
tangent spheres, then 
\beql{102}
\sum _{j=1}^{n+2} b_j^2  = \frac{1}{n} ( \sum_{j=1}^{n+2} b_j )^2.
\eeq
\end{theorem}

\noindent The case $n=3$ of this result already appears in an
1886 paper of Lachlan~\cite[p. 498]{La86} and his proof is given in
the 1916 book of Coolidge~\cite[p. 258]{Co16}.
Thus in calling this result the ``Soddy-Gossett theorem'' we are
continuing the tradition that theorems are often not named for
their first discoverers, cf. Stigler~\cite{Sti80}.
Proofs of the $n$-dimensional theorem appear in Pedoe~\cite{Pe67}
and Coxeter~\cite{Cox68b}.
Pedoe observes that this result is actually a theorem~\footnote{The
theorem depends on  the fact that the  number of
real circles, simultaneously tangent to each of $n+1$ mutually
tangent real circles with distinct tangents, is exactly two. The
total number of complex circles with this tangency property is two in
dimension $n=2$ but typically exceeds two
in dimensions $n \ge 3$.}
of real algebraic geometry,
rather than of complex algebraic geometry, in dimensions $3$ and above.

In this paper we present some very simple and elegant
extensions of these
results, which involve the centers of the circles. 
 We show that there are
relations, similar to \eqn{102}, involving the centers,
together with the curvatures, in the combination 
curvature$\times$center.
Furthermore, all these relations
generalize  to arrangements of $n+2$ mutually tangent
$(n - 1)$-spheres in $n$-dimensional Euclidean, spherical
and hyperbolic spaces, and have  a matrix formulation. 
In the process  we recover spherical and
hyperbolic analogues of the Soddy-Gossett Theorem;
these  were first obtained by Mauldon~\cite{Mau62} in 1962.
There is a vast literature
on this subject, spanning two centuries,
but (so far) we have not found our matrix
formulations in the literature.
In spirit the ideas trace back at least to Wilker \cite[p. 390]{Wi81},
see the remark at the end of \S4.  

%
%

\section{The Complex Descartes Theorem}
\setcounter{equation}{0}

\begin{figure}[htbp]
\centerline{\epsfxsize=4.0in \epsfbox{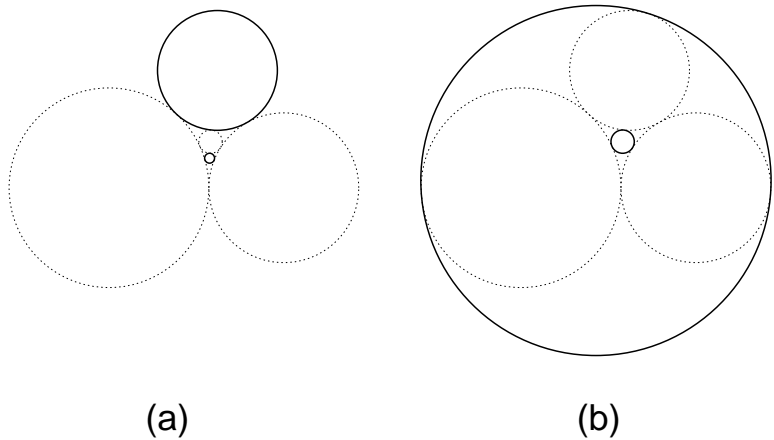}}
\caption{Circles Tangent to Three Tangent Circles}~\label{fig2}
\end{figure}

Given any three mutually tangent circles,  with curvatures 
$b_1 , b_2 , b_3$, there are
exactly two other circles that are tangent to each of these; 
each gives a  four-circle
Descartes configuration.  See Figure~\ref{fig2} for the
possible arrangements of the resulting five circles; the
three initial circles are given by dotted lines.  

The curvatures of these two new circles 
are the roots of the quadratic equation \eqn{101} 
(treating $b_4$ as the variable). 
Suppose these roots are $ b_4$ and $b_4'$.   Both can be 
positive, as in Figure 2(a), or one may be negative as in 
Figure 2(b). From \eqn{101} we have
\beql{201}
b_4 + b_4' = 2(b_1 + b_2 + b_3 ).
\eeq

\begin{figure}[htbp]
\centerline{\epsfxsize=6.0in \epsfbox{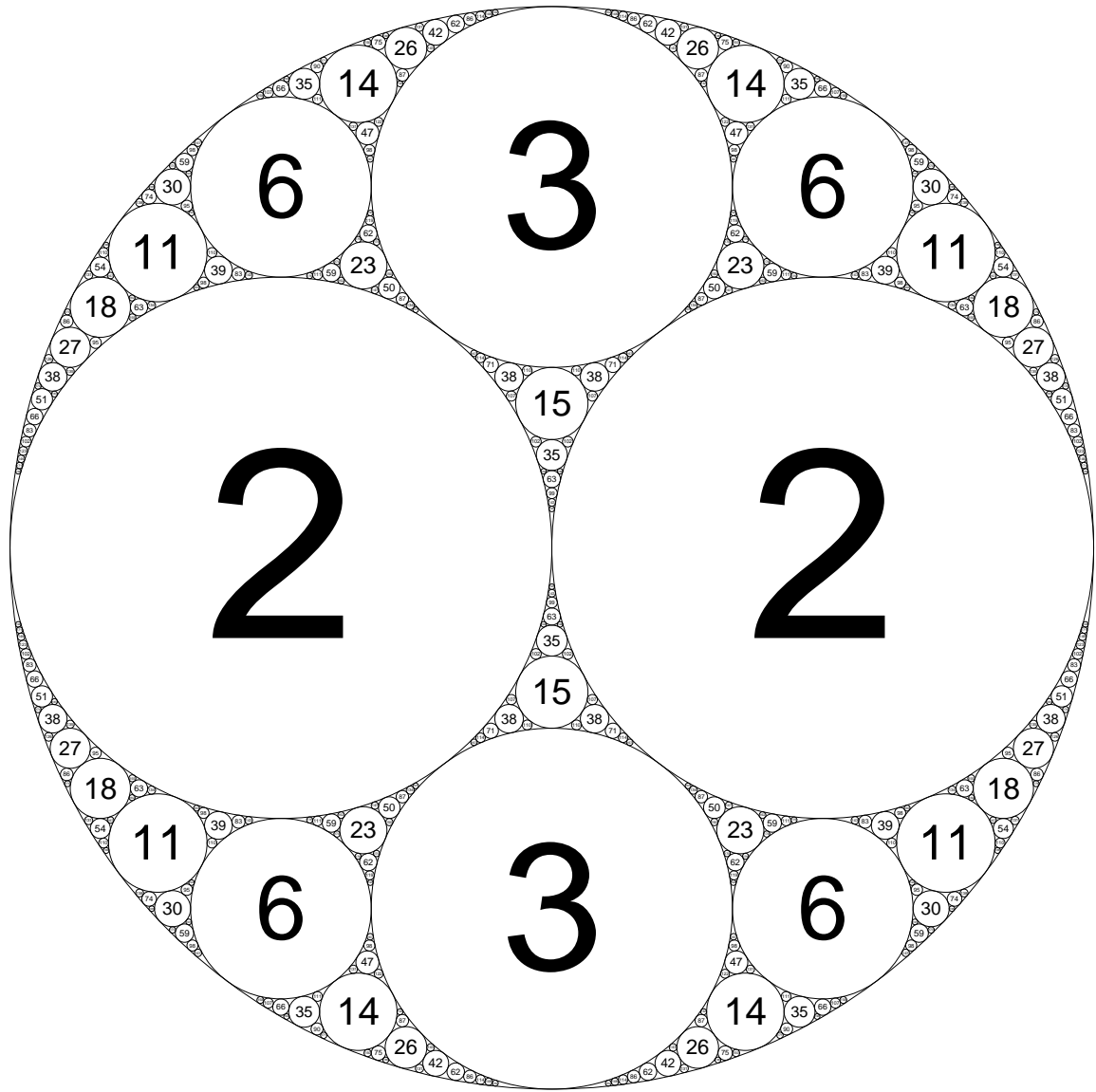}}
\caption{An Apollonian packing}\label{std}~\label{fig3}
\end{figure}

\noindent Thus, starting from a Descartes configuration, 
we can select any one of the 
four circles and replace it by the other circle that is tangent to the
remaining three; this gives a new Descartes configuration.   The new
curvature can be obtained from the original four by using \eqn{201}.  This
construction can be repeated indefinitely.  We arrive at a
packing of circles that fills either (i) a single circle, as for example 
in Figure~\ref{std},
or (ii) a strip between two parallel lines, 
or (iii) a half-plane, or (iv) the whole plane.    
Such a figure is called
an {\em Apollonian packing,} in honor of Apollonius of Perga, who considered 
(about 200 BC) the eight circles that are tangent to each of three given 
circles in general position, cf. Kasner and Supnick~\cite{KS43}.  
An Apollonian packing is completely
specified by any three mutually tangent circles in it.
 
In constructing the Apollonian
packing pictured in Figure 3, we started with four circles with vector of
curvatures $(-1,2,2,3)$.
Each circle has been labelled with its curvature; we 
notice that these are all integers.    It is clear from \eqn{201} that once 
we have a Descartes configuration with all curvatures integral, 
then in this construction all the curvatures in the packing will be integers.

In 1998, one of us, while computing Figure~\ref{fig3}, 
with the center of the outer circle located
at the origin,  noticed that the centers of 
all the circles are rational; in fact in this figure, if a circle has 
curvature $b$ and center $(x,y)$  then (it appeared) 
$b x$ and $b y$ are always 
integers.  Following this clue, we were led to the following generalization 
of \eqn{101}, in which the centers are taken to be the complex numbers 
$z_j = x_j  + iy_j$.

\begin{theorem}[Complex Descartes Theorem]~\label{th21}
Any Descartes configuration of four mutually tangent circles,
with curvatures $b_j$ and centers $z_j = x_j + i y_j$ satisfies
\beql{202}
\sum_{j=1}^4 (b_j z_j)^2 = \frac{1}{2} ( \sum_{j=1}^4 b_j z_j )^2.
\eeq
\end{theorem}
 
Notice that the relation \eqn{202} is of the same form as the original
Descartes' relation \eqn{101}.
The complex Descartes theorem implies both the Descartes circle
theorem \eqn{101} and a third relation
\beql{203}
\sum_{j=1}^4 b_j(b_j z_j) = \frac{1}{2} ( \sum_{j=1}^4 b_j)
( \sum_{j=1}^4 b_j z_j),
\eeq	
These results are  obtained by replacing $z_j$ by $z_j + w$ in \eqn{202}, 
where
$w$ is an arbitrary complex number, and identifying coefficients of
powers of $w$. 

The complex Descartes theorem also implies a relation similar to 
\eqn{201} connecting the centers of two circles,  
each of which is tangent to
each of three given mutually tangent circles, namely:
\beql{204}
b_4 z_4 + b_4' z_4' = 2( b_1 z_1 + b_2 z_2 + b_3 z_3).
\eeq
Thus in the iterative construction of an Apollonian packing
that we described above, both the curvatures
and the centers of the new circles
can be obtained by simple linear operations (followed 
by divisions). This makes it very easy to draw figures such as
Figure 3 using the computer.

The relations in the complex Descartes theorem
can be expressed in a more elegant form using the matrix

\beql{205}
\bQ_2 :=  I_4 - \frac {1}{2} \bo_4 \bo_4^T= \frac{1}{2}\left[ \begin{array}{cccc}
             1 & -1 & -1 & -1\\
            -1 &  1 & -1 & -1\\
            -1 & -1 &  1 & -1\\
            -1 & -1 & -1 &  1 \end{array} \right],
\eeq
in which  $\bo_n$ denotes a column of $n$ $1$'s , 
and $\bQ_2$  is the coefficient matrix of the {\em Descartes quadratic form}
$$ 
Q_2(x_1, x_2,x_3, x_4) :=  \bx^T \bQ_2 \bx =
(x_1^2 + x_2^2 + x_3^2 + x_4^2) - \frac {1}{2} (x_1 + x_2 + x_3 + x_4)^2.
$$
The subscript $2$ in
$\bQ_2$ refers to
the dimensionality of the space we are considering,
 
If $\bb=(b_1, b_2, b_3, b_4)^T$ denotes  the column
vector of curvatures, and $\bc = (b_1z_1, b_2z_2, b_3z_3, b_4z_4)^T$,
then 
the Descartes theorem asserts that
\beql{206}
 \bb^T \bQ_2 \bb = 0,
\eeq  
and the complex Descartes theorem asserts that
\beql{207}
 \bc^T \bQ_2 \bc = 0.
\eeq 
The complex Descartes Theorem does not completely
characterize Descartes configurations in the Euclidean plane.
There is a slightly stronger result which does, namely:

\begin{theorem}[Extended  Descartes Theorem]~\label{th22}
Given a  configuration of four oriented circles with non-zero curvatures 
$(b_1,~b_2,~b_3,~b_4)$ and 
centers $\{(x_i, y_i)~: 1 \leq i \leq 4\} $, let  
$\bM$ be the $4 \times 3$ matrix 
\beql{208}
\bM := 
\left[
\begin{array}{ccl}
b_1 & b_1x_1 & b_1y_1 \\
b_2 & b_2x_2 & b_2y_2 \\
b_3 & b_3x_3 & b_3y_3  \\
b_4 & b_4x_4 & b_4y_4
\end{array}
\right] .
\eeq
Then this  configuration 
is an oriented  Descartes configuration 
if and only if 
\beql{209}
\bM^T  \bQ_2 \bM =
     \left[ \begin{array}{ccc}
           0 & 0 & 0 \\
           0 & 2 & 0 \\
           0 & 0 & 2 \end{array} \right].
\eeq
If one or two curvatures $b_i$ are zero, and the corresponding centers 
are infinite, then  $\bM$ can be
defined in such a way that this matrix identity remains true.
\end{theorem}

The complex Descartes theorem follows from this result by
applying it to  the vector $\bc = \bx + i \by$, where $\bx$ and $\by$ are
the second and third columns of $\bM.$
The extended Descartes Theorem
gracefully generalizes to $n$-dimensions,
which we turn to next.

%
%

\section{Descartes Configurations in $n$- Dimensional Euclidean Space}
\setcounter{equation}{0}

An {\em $n$-dimensional oriented Descartes configuration}  consists of
$n+2$ mutually tangent oriented $(n - 1)$- spheres $S_i$ in $n$-dimensional
space $\rr^n,$ 
having distinct
tangencies, with the orientations compatible in the sense that
all interiors are disjoint, either with the given orientation or
with the reversal of all orientation vectors. Here we suppose that
$n \geq 2$; the one-dimensional case is treated in the
concluding section.
We often regard a hyperplane as a limiting case
of a sphere, having zero curvature, with orientation given by
a unit normal vector. In what follows an ``oriented sphere'' includes
the hyperplane case unless otherwise stated.

The Soddy-Gossett theorem \eqn{102} relates the curvatures of 
such a configuration of mutually tangent $n$-spheres, and can be written
$$
 Q_n(\bb) := \bb^T \bQ_n \bb = 0,
$$
where $\bb = (b_1 , \ldots , b_{n+2})^T$ and 
$Q_n(\bx)=\bx^T \bQ_n \bx$ is the 
{\em $n$-dimensional Descartes quadratic form}
whose associated symmetric
$(n+2) \times (n+2)$ matrix $\bQ_n$ is
\beql{300c}
\bQ_n :=  I_{n+2} - \frac {1}{n} \bo_{n+2} \bo_{n+2}^T.
\eeq
The Soddy-Gossett theorem has a converse.

\begin{theorem}[Converse to Soddy-Gosset Theorem]~\label{th30}
If $\bb=(b_1,..., b_{n+2})^T$ is a \\
nonzero real column vector
that satisfies
\beql{300a}
\bb^T \bQ_n \bb = 0,
\eeq
then there exists an oriented Descartes configuration whose
oriented curvature vector is $\bb.$ 
Furthermore any  two oriented
Descartes configurations having the same oriented curvature vector 
are congruent; that is, there is a Euclidean motion taking one
to the other.
\end{theorem}

A Euclidean motion is one that preserves angles and distances; it
includes reflections, which reverse
orientations.  We do not know an easy proof of this
result; a proof appears in \cite{GLMWY13}.

The geometry of Descartes configurations is encoded in the
curvature vector $\bb.$
If all $b_i$ are non-zero and 
$\sum_{j = 1}^{n+2} b_j  > 0$, then one of the following
holds: (i) all of $b_1, ~b_2,...,~b_{n+2}$ are positive;
(ii) $n+1$ are positive and one is negative; (iii) $n+1$ are positive 
and one is zero; or (iv) $n$ are positive and equal and the
other two are zero. These four cases correspond respectively
to the following configurations of mutually tangent spheres:
(i) $n + 1$ spheres, with another in the curvilinear simplex that
they enclose; (ii) $n+1$ spheres inscribed inside another larger sphere;
(iii) $n+1$ spheres with one hyperplane (the $(n + 2)$-nd ``sphere''),
tangent to each of them;
(iv) $n$ equal spheres with two common parallel tangent planes.

\begin{defi}~\label{de31}
{\em
Given an oriented sphere $S$ in $\rr^n$, its
{\em curvature-center coordinates} consist of the $(n+1)$-vector
$\bm(S)$ given by 
\beql{300d}
\bm(S) = ( b, b x_1, ... , b x_n)
\eeq
in which $b$ is the signed curvature of $S$ (assumed nonzero) and 
$\bx(S) = \bx  = (x_1, x_2, ..., x_n)$ is its center.
For the ``degenerate case'' of an  
oriented hyperplane $H$ 
its {\em curvature-center coordinates} $\bm(H)$ are defined to be 
\beql{300e}
\bm(S) = (0,\bh)
\eeq
where $\bh :=  (h_1, h_2, \ldots, h_n)$ is the unit normal vector
giving the orientation of the hyperplane.
}
\end{defi}

Curvature-center coordinates are not quite a global
coordinate system, because they do not always uniquely
specify an oriented sphere. Given $\bm \in \rr^{n+1}$, 
if its first coordinate 
$a \neq 0$ then there exists a unique sphere having
$\bm = \bm(S)$. But if  $a = 0$, the hyperplane case,
there is a hyperplane if and
only if $\sum h_i^2 = 1,$ and in that case there is a pencil of
hyperplanes  having the given value $\bm$, which differ from each other
by a translation.

\begin{theorem}[Euclidean Generalized Descartes Theorem]\label{th31}
Given a configuration of \\
 $n + 2$ oriented
spheres $S_1, S_2, \ldots S_{n+2}$  in $\rr^n$ (allowing
hyperplanes), let $\bM$ be 
the $(n + 2) \times (n + 1)$ matrix whose $j$-th row are the
curvature-center coordinates 
$\bm(S_j),$
of the $j$-th sphere.  
If this configuration is an oriented Descartes configuration 
then
\beql{303}
\bM^T \bQ_n  \bM  = \left[ \begin{array}{cc}
           
 	    0   &  0   \\
            0   & 2I_n    \end{array} \right] = \mbox{diag}(0,~2,~2,...,~2).
\eeq
Conversely, any real solution $\bM$ to this equation is the matrix of
a unique oriented Descartes configuration.
\end{theorem}

The curvature-center coordinate matrix $\bM$ of an oriented 
 Descartes configuration determines it
uniquely even if it contains hyperplanes, because
the other spheres in the configuration give enough information
to fix the locations of the hyperplanes.
This result contains the Soddy-Gossett theorem
as its (1,1)- coordinate. We derive the ``if'' part of this
theorem from the next result, proved in \S5. However the converse
part of this theorem seems more difficult, and we do not prove it
here. A proof appears in \cite{GLMWY13}.

We proceed to   
a further generalization, which extends the $ (n + 2) \times (n + 1)$
matrix $\bM$ to an $ (n + 2) \times (n + 2)$ matrix $\bW$ obtained by
adding an additional column. This
augmented matrix incorporates  information about
two oriented Descartes configurations, the original one and one obtained
from it by inversion in the unit sphere, as we now explain.  The
definition of $\bW$ may seem pulled out of thin air, but  
in the next two sections we will show that it
naturally arises from an analogous result
in spherical geometry, which is how we discovered it.

In $n$-dimensional Euclidean space, the operation of 
{\em  inversion in the unit sphere} replaces the point  $\bx$ by 
$\bx/{{|\bx|}^2}$, where $|\bx|^2 = \sum_{j = 1}^n x_j^2$.
Consider a general oriented sphere $S$ with center
$\bx$ and oriented radius $r$.  Then inversion in the unit sphere takes
$S$ to the sphere $\bar{S}$ with center 
$\bar{\bx} = {\bx}/(|\bx|^2 - r^2 )$
and oriented radius $\bar{r} = r/(|\bx|^2 - r^2 )$.  Note that
if ${|\bx|}^2 > r^2$, $\bar{S}$ has the same  orientation  
as $S$.   In all cases,
\beql{303f} 
\frac {\bx}{r} = \frac{\bar{\bx}}{\bar{r}}
\eeq 
and 
\beql{306aa}
\bar{b} = \frac{|\bx|^2}{r}- r.
\eeq

\begin{defi}~\label{de32} 
{\em 
Given an oriented sphere $S$ in $\rr^n$, its
{\em augmented curvature-center coordinates} are the $(n+2)$-vector
 
\beql{300g}
\bw(S) := ( \bar{b}, b, b x_1, ... , b x_n) = (\bar{b}, \bm),
\eeq
in which $\bar{b}= b(\bar{S})$,
 is the curvature of the sphere or hyperplane $\bar{S}$
obtained by inversion of $S$ in the unit sphere,
and $\bm$ are its curvature-center coordinates.  For hyperplanes we
define
\beql{300f}
\bw(H) := (\bar{b}, 0, h_1, ... , h_n) = (\bar{b}, \bm),
\eeq
where $\bar{b}$ is the oriented curvature of the sphere 
or hyperplane $\bar{H}$ obtained by inversion of $H$ in the unit
sphere.
}
\end{defi}

Augmented curvature-center coordinates provide a global
coordinate system: no two distinct oriented spheres have the
same coordinates. The only case to resolve is when $S$ is a hyperplane,
i.e.  $b=0$. The relation \eqn{303f} shows that
$(\bar{b}, b x_1, ... , b x_n)$ are the curvature-center coordinates
of $\bar{S}$, and  if $\bar{b} \neq 0$, this uniquely determines $\bar{S}$,
and then, by inversion in the unit circle, $S$. In the
remaining case $b = \bar{b} =0$ then $S = \bar{S}$ is the
unique hyperplane passing through the
origin whose  unit normal is given by the remaining coordinates.

Given a collection  $(S_1, S_2, ... , S_{n + 2})$
of $n + 2$ oriented spheres (possibly hyperplanes)
in $\rr^n$, 
the {\em augmented  matrix} $\bW$ associated with it is
the $(n+2) \times (n + 2)$ matrix whose $j$-th row is given by the
augmented curvature-center coordinates 
$\bw(S_j)$ of the $j$-th sphere.

The action of inversion in the unit sphere has a
particularly simple interpretation in augmented matrix coordinates. If $\bW$
is the augmented matrix associated to a Descartes configuration, and if
$\bW^\prime$ is the augmented matrix associated to its inversion in
the unit sphere, then
\beql{306a}
\bW = \bW^\prime \left[ \begin{array}{ccl}
           0 &   1   & 0   \\
	   1 &   0   & 0   \\
           0 &   0   & I_n    \end{array} \right],
\eeq
a result which follows from \eqn{303f}.

\begin{theorem}
[Augmented Euclidean Descartes Theorem]~\label{th32}
An oriented Descartes configuration of $n + 2$  spheres 
$\{S_i:  1 \leq i \leq n+2\}$ in $\rr^n$
has an augmented matrix $\bW$ which satisfies
\beql{307}
\bW^T \bQ_n \bW  =  \left[ \begin{array}{ccl}
           0 &  -4   & 0   \\
	  -4 &   0   & 0   \\
           0 &   0   & 2I_n    \end{array} \right].
\eeq
Conversely, any real solution $\bW$ to this matrix equation 
is the augmented matrix
of a unique oriented Descartes configuration.
\end{theorem}

The augmented Euclidean Descartes theorem
includes as special cases the ``if'' direction
of each of the theorems stated so far,
and represents our final stage of generalization
of the Descartes circle theorem in Euclidean space.
In particular, the ``if'' part of the 
Euclidean generalized Descartes theorem is just \eqn{307} with the first
row and column deleted.
In the converse direction~\footnote
{In the converse direction the Augmented Euclidean
Descartes theorem is not as strong as the converse in the
Euclidean generalized Descartes theorem, nor does it imply
the converse to the Soddy-Gossett theorem;
these results require separate proofs.}
this theorem gives a
parametrization of  all oriented Descartes configurations,
and it is  a ``moduli space'' for such configurations given as
an affine real-algebraic variety.

We discovered the augmented Euclidean Descartes theorem in
studying analogues of the Descartes theorem in non-Euclidean
geometries.  In the next section we formulate and prove such
an analogue in spherical geometry; then in \S5 we deduce 
the augmented Euclidean Descartes theorem from it.

%
%

\section{Spherical Geometry}
\setcounter{equation}{0}

The standard model for spherical geometry $\sss^n$ is the unit $n$-sphere
$S^n$ embedded in $\rr^{n+1}$ as the surface 
\beql{401}
S^n := \{y: y_0^2 + y_1^2 + \ldots +y_n^2 = 1 \} 
\eeq
with the Riemannian metric induced from the Euclidean metric 
in $\rr^{n + 1}$ by restriction.   In this model,  the distance between
two points of $\sss^n$ is simply the angle $\alpha$ between 
the radii that join the origin of $\rr^{n+1}$ to the representatives
of these points on $S^n$.   This distance $\alpha$
always satisfies $0 \leq \alpha \leq \pi$.

A {\em sphere} $C$ in this geometry 
 is the locus of points 
equidistant (at distance $\alpha$ say) from a point in $\sss^n$ called 
its {\em center.}  The quantity $\alpha = \alpha(C)$ is the
{\em spherical radius} or {\em angular radius} of $C$;
it is the angle at the origin ${\bf 0}$ of $ \rr^{n+1}$   
between a ray from
 ${\bf 0}$ to the center of $C$ and a ray from
${\bf 0}$ to any point of $C$.
 Note that there are
two choices for the center (and the angular radius)
of a given sphere; these two choices form
a pair of antipodal points of $S^n$. The choice of a center amounts to
orienting the sphere.
In this model the interior of a sphere is a  spherical cap, 
cut off by  the
intersection of the sphere $S^n$ with a hyperplane in $\rr^{n + 1}$,
so (by abuse of language) we will also 
call an oriented sphere a {\em spherical cap}.
 
The two spherical caps determined by a given sphere are called
{\em complementary} and the sum of their angular radii is $\pi.$
We define the {\em interior} of an oriented sphere to contain all points
of $S^n$ on the same side of the hyperplane 
as the center of the sphere.   If we describe a hyperplane by a linear 
form
\beql{402b}
\bF(y) = \sum _{i = 0}^n f_i y_i - f,
\eeq
normalized by the requirement 
$$ \sum_{i = 0}^n f_i^2 = 1, $$
this provides an orientation by defining a positive half-space $F(y)  > 0$.
The sphere has center $\bff := (f_0, f_1, \ldots ,f_n)$ and 
has positive radius 
if and only if $|f| < 1$.  The 
radius $\alpha$ satisfies $\cos \alpha = f$, and the interior
of the spherical cap it determines is the region where the
linear form is positive. 
A spherical cap can 
be specified either by a pair $(\bff, \alpha )$ or by the pair 
$(-\bff, \alpha - \pi )$, while $(-\bff, \pi - \alpha)$ determines
the complementary spherical cap.

A {\em spherical Descartes configuration} consists of 
$n+2$ mutually tangent spherical caps
 on the surface of the unit n-sphere, such that either (i)
the interiors of all spherical caps are mutually disjoint, or 
(ii) the interiors
of all complementary spherical caps are mutually  disjoint.

\begin{theorem}[Spherical Soddy-Gossett Theorem]~\label{th71a}
Given a spherical Descartes configuration of
 $n + 2$ mutually tangent
spherical caps $C_i$ on the $n$-dimensional
unit sphere $S^n$ embedded in $\rr^{n+1},$ with spherical radius
$\alpha_j$ subtended by the $j$-th cap, then these
spherical radii satisfy the relation
\beql{N702}
\sum_{i=1}^{n+2}(\cot \alpha_i)^2 = 
\frac {1}{n} (\sum_{i=1}^{n+2} \cot \alpha_i )^2 
 -2.
\eeq
\end{theorem}

This theorem was found by Mauldon~\cite[Theorem 4]{Mau62} 1n 1962,
as part of a more general result allowing non-tangent spheres.
He also established 
a converse: to each real solution
of \eqn{N702} there corresponds some spherical Descartes
configuration, and two spherical Descartes configuration with
the same data in \eqn{N702} are congruent configurations
in spherical geometry. 

The spherical Soddy-Gossett theorem
 is intrinsic, i.e. it depends only on the Riemannian metric
for spherical geometry, and not on the coordinate system used to
describe the manifold. 
However we shall establish  it as a special case of
a result that does depend on a particular choice of coordinate
system. If $C$ is a spherical cap with center 
$\by = (y_{0}, y_{1}, y_{2}, ... y_{n+1})$,
and angular radius $\alpha$, we define its
{\em spherical curvature-center coordinates} $\bw_+(C)$ to be the
 row vector
\beql{404a}
\bw_+(C) := ( \cot \alpha, \frac{y_0}{\sin \alpha}, \frac{y_1}{\sin \alpha},
\ldots,\frac{y_n}{\sin \alpha}).
\eeq
No two spherical caps have the same coordinates $\bw_+,$
 since $\alpha$ is uniquely determined by the
first coordinate, and then the $y_j$ are uniquely determined using
the other coordinates.

To any configuration of $n+2$ caps $C_1, \dots , C_{n+2}$
we associate the $(n+2) \times (n+2)$
{\em spherical curvature-center coordinate matrix}
 $\bW_+$ whose  $j$th row is $\bw_+(C_j)$.  

\begin{theorem}[ Spherical Generalized Descartes Theorem]~\label{th71}
Given a  configuration of $n+2$ oriented spherical caps $C_j$
that is a spherical Descartes configuration, then the
$(n+2) \times (n+2)$ matrix $\bW_+$ whose $j$-th row is  the
spherical curvature-center coordinates of $C_j$ satisfies 
\beql{402}
           \bW_+^T \bQ_n \bW_+  = 
      \left[ \begin{array}{ccc}
           -2 & 0 & 0 \\
           0 & 2 & 0 \\
           0 & 0 & 2I_n \end{array} \right] = 
 \mbox{diag}(-2, 2, 2, ... , 2).
\eeq
Conversely, any real matrix $\bW_{+}$ that satisfies this equation
is the spherical curvature-center coordinate
matrix of some spherical Descartes configuration.
\end{theorem}

The $(1,1)$-entry of the matrix relation \eqn{402} is
the spherical Soddy-Gossett theorem above.

This theorem has a remarkably simple proof, which is
based  on two preliminary
lemmas. 
Let $\bJ_n$ be the $(n+2) \times (n+2)$ matrix 
\beql{403aa}
\bJ_n =  \left[ \begin{array}{ccc}
           -1 &   0   & 0   \\
	    0 &   1   & 0   \\
            0 &   0   &  I_n    \end{array} \right]= \mbox{diag}(-1,1,...,1).
\eeq
 
\begin{lemma}~\label{le41}
(i) For any $(n + 2)$-vector $\bw_+$, there is a spherical cap 
$C$ with $\bw_+(C) = \bw_+$ if and only if
\beql{405a}
\bw_+ \bJ_n \bw_+^T = 1.
\eeq

(ii) The spherical caps $C$ and $C'$ are externally tangent if and only if
\beql{406a}
\bw_+(C) \bJ_n \bw_+(C')^T = -1.
\eeq
\end{lemma}

\pf.
(i). If $\bw_{+}$ comes from a spherical cap with center $\by$
and angular radius $\alpha$, then 
$$ \bw_+ \bJ_n \bw_+^T = \frac{ - (\cos \alpha)^2 + 
\sum_{j = 0}^n y_j^2}{(\sin \alpha)^2} = \frac { 1 - (\cos \alpha)^2}{(\sin \alpha)^2} = 1$$
so \eqn{405a} holds.

Conversely, if \eqn{405a} holds,
then one recovers a unique $\alpha$ with $0 < \alpha < \pi$
by setting $\cot \alpha := (\bw_+)_1,$ and  
one then  defines a vector $\by= (y_0, ..., y_{n+1})$ via
$y_j := \frac{(\bw_+)_j}{\sin \alpha},$ noting that $\sin \alpha \neq 0.$
The equation \eqn{405a} now implies that  ${|\by|}^2 = 1,$ 
so $\by$  lies on the unit sphere, and we have 
determined a spherical cap giving the vector $\bw_+$.

(ii). Two spherical caps with centers $\by, \by^\prime$ 
with angular radii $\alpha,
\alpha^\prime$ are externally tangent if and only if the angle between their
centers, viewed from the origin in $\rr^{n+1}$ is $\alpha + \alpha^\prime.$
Since $\by$ and $\by^\prime$ are unit vectors, this
holds if and only if 
$$
\by (\by^\prime)^T = \cos ( \alpha + {\alpha}^\prime),
$$
Now
$$\bw_+(C) \bJ_n \bw_+(C')^T = \frac{1}{\sin \alpha \sin \alpha^\prime}
(-\cos \alpha \cos \alpha^\prime + \by (\by^\prime)^T ) $$
and this  gives \eqn{406a}, using 
$\cos(\alpha + \alpha^\prime) = \cos \alpha \cos \alpha^\prime -
\sin \alpha \sin \alpha^\prime.$
~~~$\Box$

\begin{lemma}~\label{le42}
If $\bA$, $\bB$ are symmetric non-singular $n \times n$ matrices 
and $\bW \bA {\bW}^T = \bB$, then
${\bW}^T {\bB}^{-1} \bW = {\bA}^{-1}$.
\end{lemma}

\pf.
The matrix $\bW$ is non- singular since $\bB$ is
non-singular. Now invert 
both sides, then  multiply  on the left 
by ${\bW}^{T}$ and on the right by ${\bW}$.
~~~$\Box$

\noindent\paragraph{Proof of the Spherical Generalized Descartes Theorem.}
If the caps $C_j$ touch externally, we have from Lemma \ref{le41} that
\beql{404}
\bW_+ \bJ_n \bW_+^T = 2\bI_{n+2} - \bo_{n+2} \bo_{n+2}^T = 2\bQ_n^{-1}.
\eeq
Then applying  Lemma \ref{le42} (with $\bA = \bJ_n$ and $\bW = \bW_+$) 
we obtain 
\beql{405}
\bW_+^T \bQ_n \bW_+ = 2{\bJ_n}^{-1} = 2\bJ_n.
\eeq

Conversely, \eqn{405} implies \eqn{404}, by Lemma 4.4. Looking
at the  diagonal elements of $\bW_+ \bJ_n \bW_+^T,$ 
which are all  $1$'s, 
Lemma~\ref{le41}(i) guarantees that
the $j$-th row of $\bW_{+}$ is a vector $\bw_+(C_j)$ for some
(uniquely determined) spherical cap $C_j$, and the 
off-diagonal elements all being $-1$ shows by Lemma~\ref{le41}(ii) that
the caps touch externally pairwise, so form a spherical Descartes 
configuration.
~~~$\Box$

\noindent\paragraph{Remark.} Wilker~\cite[pp. 388-390]{Wi81}
came tantalizingly close to obtaining the spherical generalized
Descartes theorem. He termed a spherical
Descartes configuration a ``cluster'',
and introduced spherical curvature-center coordinates.
In a remark he noted our Lemma~\ref{le41} and stated equation \eqn{404}.
However he did not invert his formula, via Lemma~\ref{le42}
and so failed to formulate a result in terms of the Descartes
quadratic form. 

%
%

\section{Stereographic Projection and the Augmented Euclidean
Descartes Theorem}
\setcounter{equation}{0}

We derive  the augmented Euclidean Descartes  theorem
from the spherical Generalized Descartes
theorem, using stereographic projection. 
The resulting
derivation is reversible, so the spherical Generalized
Descartes theorem and the augmented Euclidean Descartes
theorem may be viewed as equivalent results.

\begin{figure}[htbp]
\centerline{\epsfxsize=4.0in \epsfbox{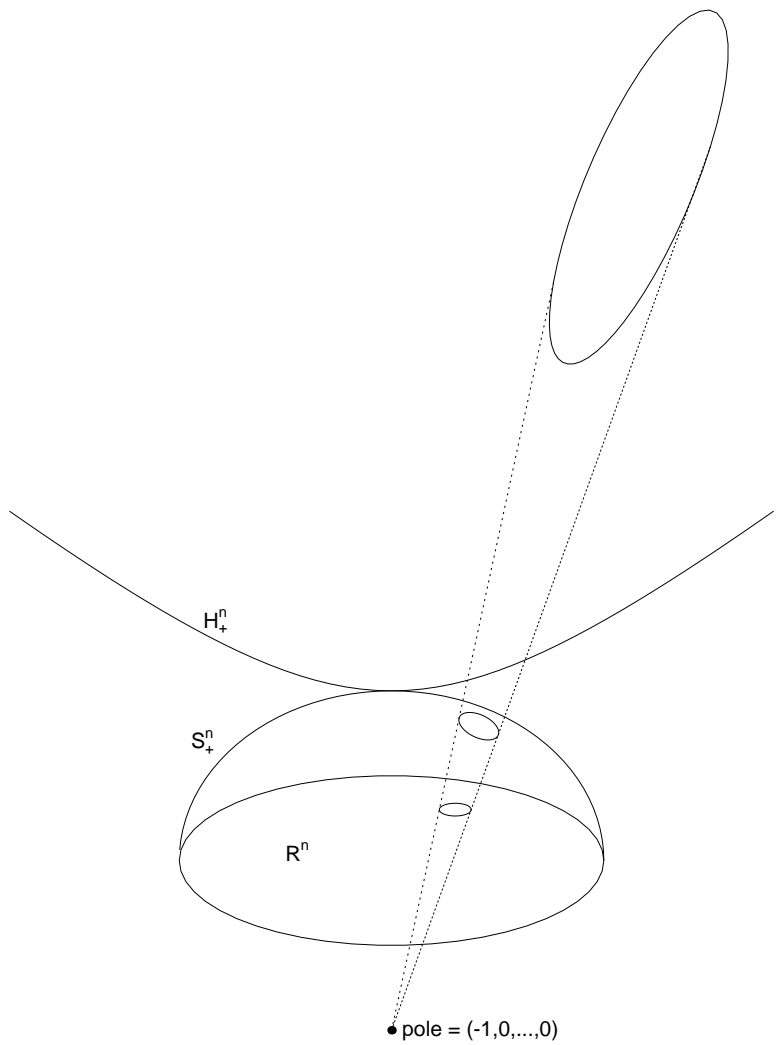}}
\caption{Stereographic projection-hyperplane,
 sphere and hyperboloid}\label{fig4}
\end{figure}

Consider the unit sphere in $\rr^{n+1}$, 
given by $\sum_{i=0}^{n} y_i^2 = 1$.  
Points on 
this sphere can be mapped into the plane $y_0 = 0$ by 
stereographic projection 
from the ``south pole''
$(-1,0,...,0)$, see Figure~\ref{fig4}. (The hyperboloid
in the figure will be used later.)

 This mapping $(y_0,...,y_n) \to (x_1,...,x_n)$ is given by
$$ 
  x_j = \frac{y_j}{1+y_0},~~~ 1 \le j \le n.
$$ 

The spherical cap $C$ with center $(p_0,...,p_n)$ and angular
radius  $\alpha$ is the intersection of the unit sphere with the plane
$$
    \sum_{j = 1}^n  p_j y_j = \cos{\alpha}.
$$
The sterographic projection of this cap in 
the hyperplane $y_0 = 0$ is the (Euclidean) 
sphere $S$ with 
center $(x_1,...,x_n)$ and radius $r$, where
$$
   x_j = \frac{p_j}{p_0 + \cos{\alpha}},~~~~ 1 \leq j \leq n,
$$
and  
$$  
r = \frac{ \sin{\alpha}}{p_0 + \cos{\alpha}}.
$$
If the boundary of the cap $C$ contains the south pole,  the corresponding
sphere $S$ has infinite radius, i.e. it is a hyperplane.

\noindent\paragraph{Proof of the Augmented Euclidean Descartes Theorem.}
The spherical coordinates of the spherical cap $C$ are given by
the row-vector
$$
\bw_+(C) = ( \cot \alpha, \frac{p_0}{\sin \alpha}, \frac{p_1}{\sin \alpha},
\ldots,\frac{p_n}{\sin \alpha}),
$$
We relate this to the augmented Euclidean coordinate
vector $\bw(S)$ associated with the corresponding
projected sphere $S$ in the plane $y_0 = 0$, given by \eqn{300g}.  We have 
 ${x_j}/{r} = {p_j}/{\sin \alpha}$,  
$b = 1/r = \cot \alpha +  \frac {p_0}{\sin \alpha}$, and we find
$$
\bar{b} = \cot \alpha - \frac {p_0}{\sin \alpha}.
$$
Thus  
$$\bw(S) = (\cot\alpha - \frac{p_0}{\sin \alpha}, 
\cot \alpha + \frac{p_0}{\sin \alpha}, \frac{p_1}{\sin \alpha},...,
\frac{p_n}{\sin \alpha}) =\bw_+(C)\bG,$$ 
where
\beql{412a}
\bG = \left[ \begin{array}{ccc}
           1 &  1 & 0 \\
          -1 &  1 & 0 \\
           0 &  0 & I_n \end{array} \right].
\eeq

Suppose we have a configuration of $n+2$ spherical caps 
$C_1, ... , C_{n+2}$ on the unit sphere.
These stereographically project into a configuration of 
Euclidean spheres 
$S_1, ... , S_{n+2}$ 
in the equatorial plane $y_0=0$, and conversely every 
configuration of Euclidean spheres lifts to a configuration
of spherical caps.
The map sends spherical Descartes configurations to
Euclidean Descartes configurations.
   We assemble the corresponding rows $\bw_+(C_j)$,
$\bw(S_j)$ into matrices $\bW_+$ and $\bW$,
respectively.  Then 
\beql{413}
\bW = \bW_+ \bG,
\eeq 
and, using the Spherical Generalized
Descartes Theorem~\ref{th71}, we have
$$
\bW^T \bQ_n \bW =  \bG^T \bW_+^T \bQ_n  \bW_+ \bG = 
\bG^T \mbox{diag}(-2,2,...,2) \bG =  \left[ \begin{array}{ccl}
           0 &  -4   & 0   \\
	  -4 &   0   & 0   \\
           0 &   0   & 2I_n    \end{array} \right],
$$
which proves the augmented Euclidean Descartes Theorem~\ref{th32}.~~~$\Box$

%
%

\section{Hyperbolic Geometry}
\setcounter{equation}{0}

There are many models
of hyperbolic space $\hh^n$, of which the three
most common are the 
(Poincar\'{e}) unit ball model, the  half space model, 
and the hyperboloid model.
(In two \\
dimensions we say ``unit disk'' and ``half-plane'' for the 
first two models.)
The unit ball and half-space models are described in many places, e.g.
Beardon~\cite{Be83} and Berger~\cite[Chapter 19]{Ber87}. 
The hyperboloid model, which is less well known, but which is in some 
ways simpler than the others,
is described in Beardon~\cite[Section 3.7]{Be83},
Reynolds~\cite{Re93} and Ryan~\cite{Ry86}.
The unit ball and half-space models
are embedded in  $\rr^n$,
though with different metrics, while the hyperboloid model
is embedded in $\rr^{n+1}$, endowed with a Minkowski metric.  
Here we
need only the unit ball and hyperboloid models. 
A {\em  sphere} in hyperbolic $n$-space $\hh^n$ is defined as the locus of
points equidistant (in hyperbolic metric) from some fixed point in $\hh^n$,
the center.

The unit ball model
consists of the points  $(y_1, \ldots, y_n)$ in $\rr^n$
with $\sum_{j = 1}^n y_i^2 < 1,$  with
the ideal boundary being $\sum_{j = 1}^n y_j^2 =1.$
In this model, the hyperbolic metric is
$$
ds^2 = (dy_1^2 + \dots + dy_n^2)/(1- \sum_{j = 1}^n y_j^2)^2. 
$$
and the hyperbolic distance between two points
$\by, \by'$ satisfies
\beql{N729}
\cosh (d (\by, \by')) = 
\left( (1+\sum_{j = 1}^n y_j^2)(1+ \sum_{j = 1}^n {y_j'}^2) -
4 \sum_{j = 1}^n y_j y_j'  \right) / \left( (1- \sum_{j = 1}^n y_j^2)
(1- \sum_{j = 1}^n {y_j'}^2) \right).
\eeq
In this model  
a hyperbolic sphere (of finite radius) is a Euclidean
sphere contained strictly inside the unit ball; however its hyperbolic
center and hyperbolic radius usually differ from the Euclidean ones.

Points in the hyperboloid model are represented in $\rr^{n + 1}$
as points on the upper sheet $H_+^n$ given by
$u_0 > 0$ of the two sheeted-hyperboloid 
$H_{\pm}^n$     
cut out by the equation 
$$ u_0^2 = 1 + u_1^2 + \cdots + u_n^2,$$
where $H_{\pm}^n= H_+^n \cup H_-^n$ with $H_-^n = - H_+^n.$
However the Riemannian metric is not that induced from
the Euclidean metric on $\rr^{n + 1},$ but rather is
induced from the Minkowski metric
$$ds^2 = - du_0^2 + du_1^2 + ... + du_n^2$$
on this space, cf. Beardon~\cite[p. 49]{Be83}.
The formula for the {\em hyperbolic distance} $d( \bu, \bu')$ in
this metric is given by
\beql{N730}
\cosh (d( \bu, \bu')) = u_0 u_0' - u_1 u_1' - \ldots - u_n u_n',
\eeq
see Reynolds~\cite[formula (6.10)]{Re93}.
One can go between
the hyperboloid model and the ball model by the change
of variables
$$y_j = \frac {u_j}{1 + u_0}, \qquad\mbox{for}~~1 \leq j \leq n,$$
and in the opposite direction by 
$$
u_0 = \frac{2}{\Delta} - 1 \qquad\mbox{and}~~  
u_j = \frac {2 y_j }{ \Delta}, ~~ 1 \le j \le n, 
$$
where 
$$\Delta = 1-\sum_{j = 1}^n u_j^2.$$

From  \eqn{N730} we see that in this model a hyperbolic  sphere is
represented by the intersection of $H_+^n$ with a hyperplane
$\bG(u)=0$, where
\beql{N730b}
\bG(u) =  g_0u_0 - \sum _{i = 1}^n g_i u_i - g,
\eeq
where $g > 1$ and we require $\bG$ to be normalized by the requirement 
that 
\beql{N730c}
 g_0^2 = 1 + \sum_{i = 1}^n g_i^2,
\eeq
i.e. the point $\bg := (g_0, g_1, \ldots , g_n)$ lies on $H_+^n$.  
This intersection is typically a (Euclidean)
 ellipsoid.  The center of the sphere is
$\bg$,  and its radius $d$ satisfies $\cosh d = g$.   As in the spherical 
case, we define the interior of the hyperbolic
sphere to be the region on the same 
side of the plane $\bG(u)=0$ as the center.  

If we consider a general hyperplane $\bG(u)=0$ which has the  normalized
form  \eqn{N730c}, intersecting 
the complete hyperboloid 
$H_+^n \cup H_-^n$, the intersection may be
empty, a single point, an ellipsoid (or sphere) on either sheet, a paraboloid
on either sheet, or a two-sheeted hyperboloid.  For later use we term 
these possibilities 
{\em virtual hyperbolic spheres}, except for  the empty set or a point.
We also  assign them an {\em orientation} given by the sign of the
constant term $g$ in the  
associated
linear form $\bG(u)$. (Replacing $\bG(u)$ by $-\bG(u)$ reverses
the orientation.) Only
the points in $H_+^n$ correspond to real points in $\hh^n$.  
In the two-dimensional case parabolas on the upper sheet correspond to 
{\em horocycles}; in the unit disc model these are (Euclidean) circles that 
are tangent to the bounding circle.    They have infinite radius.  In the 
disc model their centers are on the bounding circle, and in the hyperboloid 
model their centers are at infinity.  The boundary of the disc model 
corresponds  to a circle at infinity in the hyperboloid model.
 
An {\em oriented hyperbolic Descartes configuration} is any
set of $n+2$ mutually tangent oriented hyperbolic $(n - 1)-spheres$ 
in $\hh^n$, having the property that 
either (i) all interiors of the spheres are
disjoint, or (ii) the interiors of each pair of spheres intersect in
a nonempty open set.
In the following result we also allow (oriented) Descartes configurations
which include those virtual hyperbolic spheres of infinite
radius which in the ball model correspond to
Euclidean spheres lying entirely inside the closed ball and tangent to 
its boundary. The ideal boundary itself 
forms a single limiting $(n - 1)$-sphere in
this sense.

\begin{theorem}[Hyperbolic Soddy-Gossett Theorem]~\label{th72a}
The oriented hyperbolic radii \\
$\{s_j~:~ 1 \le j \le n+2 \}$
of an oriented Descartes configuration of
$n + 2$ spheres in hyperbolic space $\hh^n$
satisfy the relation
\beql{711a}
\sum_{j = 1}^{n + 2} (\coth{s_j})^2 = \frac{1}{n} 
(\sum_{j = 1}^{n+2} \coth{s_j})^2 + 2.
\eeq
\end{theorem}

This result was found by Mauldon~\cite{Mau62}.
The hyperbolic Soddy-Gossett theorem is intrinsic, 
depending only on the hyperbolic
metric. We derive it as a special case of a result which does depend
on a specific coordinate system, namely that given above for the
the hyperboloid model.

If $S$ is a hyperbolic sphere in $H_+^n$  with center 
$\bu = (u_{0}, u_{1}, u_{2}, ..., u_{n+1})$,
and hyperbolic  radius $s_j$, we define its
{\em hyperbolic curvature-center  coordinates} $\bw_-(S)$ to be the
 row vector
\beql{504a}
\bw_-(S) :=( ~\coth s  , \frac{u_0  }{\sinh s} , \frac{u_1}{\sinh s}, 
..., \frac{u_n}{\sinh s}).
\eeq
To a configuration of $n+2$ hyperbolic spheres  $S_1, \dots , S_{n+2}$
we associate the $(n+2) \times (n+2)$
matrix $\bW_{-}$ whose  $j$th row is $\bw_{-}(S_j).$

\begin{theorem}[Hyperbolic Generalized Descartes Theorem]~\label{th72}
Given a  configuration of \\
$(n+2)$ oriented hyperbolic spheres which
is a hyperbolic Descartes configuration, then 
the  associated matrix $\bW_{-}$ whose rows are the 
hyperbolic curvature-center coordinates of the spheres
satisfies
\beql{N708}
       \bW_{-}^T \bQ_n \bW_{-} = 
       \left[ \begin{array}{ccl}
           2 &   0   & 0   \\
	   0 &  -2   & 0   \\
           0 &   0   & 2I_n    \end{array} \right]= 
 \mbox{diag}(2, -2, 2, ... , 2) .
\eeq
\end{theorem}

The converse of Theorem~\ref{th72} does not hold,
because some matrices
$\bW_{-}$ satisfying \eqn{N708} do not correspond to hyperbolic
Descartes configurations. We can obtain a converse by  
allowing ``virtual Descartes  configurations'' 
that lie on both sheets of the
hyperboloid. One simply defines a {\em virtual Descartes configuration}
to be the image on the two-sheeted hyperboloid resulting from stereographic
projection through $(-1, 0, 0,...,0)$
of any Descartes configuration on the unit sphere. The resulting 
hyperbolic coordinate matrix
$\bW_{-}$ is to be defined by \eqn{709a}. 
One can define a  ``virtual (oriented) hyperbolic sphere'', and
define its  oriented radius and center   
using the formulas following \eqn{N730b}; 
the center and oriented radius of some ``virtual hyperbolic spheres''  may
then be (non-real) complex numbers, although by definition their
hyperbolic curvature-center coordinates will be real.

Theorem~\ref{th72} is readily deducible from the 
spherical generalized Descartes theorem via
stereographic projection through the ``south pole''
 $(-1, 0, 0,...0)$ in $\rr^{n+1},$
mapping oriented Descartes configurations on the upper sheet ${H_+}^n$
of the hyperboloid
to spherical Descartes configurations which lie entirely
on the upper hemisphere${S_+}^n$ of the unit sphere. 
See Figure ~\ref{fig4}
in \S5.

In applying stereographic projection, the locus of a hyperbolic
$(n-1)$-sphere on the hyperboloid is mapped to the locus of a
spherical $(n-1)$-sphere on the unit $n$-sphere, and also to the
locus of a Euclidean $(n-1)$-sphere in the plane $x_0 = 0.$
Note however that the hyperbolic center, the spherical center of the 
associated spherical
cap and the Euclidean center of the Euclidean
sphere are typically all distinct in the sense that they
usually lie on three different lines through the ``south pole''
$(-1,0, ...,0)$ in $\rr^{n+1}.$

One can show that, if $\bW_{-}$ is defined as above, and if
$\bW_{+}$ are the spherical Descartes coordinates associated to
the resulting Descartes configuration on the sphere as in \S4,  then they are
related by
\beql{709a}
\bW_{-} = \bW_{+}\left[ \begin{array}{ccl}
           0 &   1   & 0   \\
	   1 &   0   & 0   \\
           0 &   0   & I_n    \end{array} \right].
\eeq
We omit a proof of this fact, which can be carried out along the lines of \S5.
Given it, one immediately deduces Theorem~\ref{th72} from Theorem~\ref{th71},
plus a converse if ``virtual Descartes configurations'' are included.

%
%

\section{ Apollonian Packings}
\setcounter{equation}{0}

Using stereographic projection we have a recipe to pass between
Euclidean, spherical and hyperbolic Descartes configurations. 
It gives a one-to-one correspondence between configurations
$\bW, \bW_{+},$ and $\bW_{-}$ given by \eqn{412a} and \eqn{709a},
namely
\beql{801}
\bW = \bW_{+} \left[ \begin{array}{ccc}
           1 &  1 & 0 \\
           -1 &  1 & 0 \\
           0 & 0 & I_n \end{array} \right] = 
        \bW_{-}\left[ \begin{array}{ccc}
           -1 &  1 & 0 \\
            1 &  1 & 0 \\
            0 &  0 & I_n \end{array} \right].
\eeq
Here we have extended the definition of Descartes configuration
to ``virtual Descartes configurations''
in the hyperbolic case to be configurations on both sheets of
the two-sheeted hyperboloid. This recipe clearly lifts to
Apollonian packings. 

Since the spherical and hyperbolic Soddy-Gossett theorems
involve quadratic forms, an analogue of the relation
\eqn{201} holds in spherical and hyperbolic
geometry, permitting the easy calculation of the
``curvatures'' $\cot \alpha$ (resp. $\coth s$) of circles
in spherical (resp. hyperbolic) packings. 
There is a notion of 
``integral Apollonian circle packing'' for such 
``curvatures'' which makes sense in 
spherical and hyperbolic geometry. Furthermore, analogues of
the relation \eqn{204} hold in spherical and hyperbolic
geometry as well, permitting the easy calculation of the
centers
in spherical (resp. hyperbolic) Apollonian packings. 

\begin{figure}[htbp]
\centerline{\epsfxsize=6.0in \epsfbox{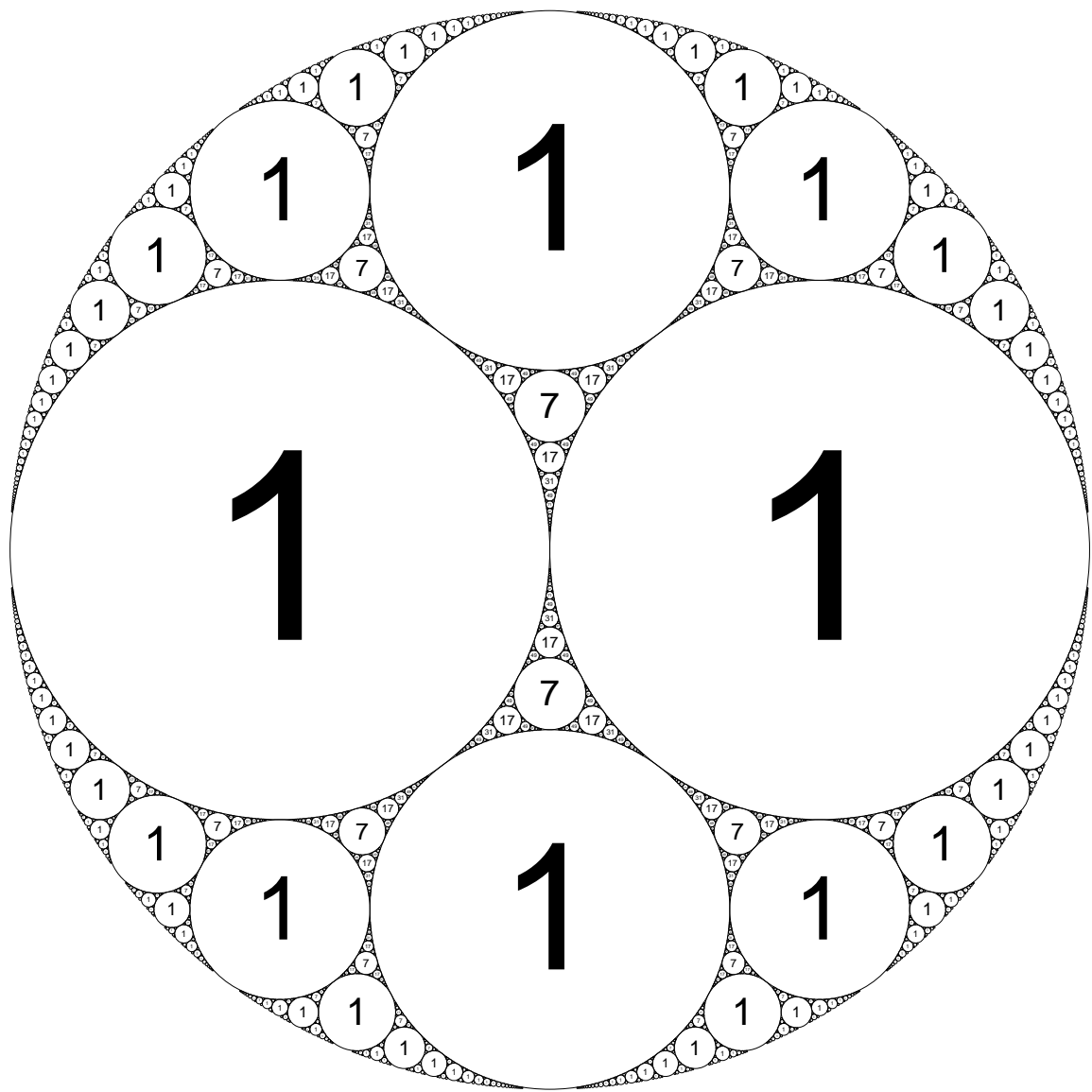}}
\caption{A hyperbolic Apollonian packing}\label{fig5}
\end{figure}

Thus the standard Euclidean Apollonian  packing 
pictured  in Figure~\ref{fig3}, with center at the origin, has
a corresponding  hyperbolic packing obtained by stereographic projection
in which the bounding outer circle
in the packing is the ``absolute'' in the unit ball model of the
hyperbolic plane, and in which 
the $\coth(r)$'s are all integers,
but {\em not} the same integers as in the Euclidean packing,
calculated using $\bW_{-}$ in \eqn{801}.
See Figure~\ref{fig5}.

\noindent Those circles that are tangent to the bounding
 circle are known as horocycles, and have 
infinite hyperbolic radius, so the corresponding value of $\coth(r)$ is
$1.$ This explains the large number of circles
assigned the value $1$ in Figure~\ref{fig5}, namely all those that touch
the outer circle. 

\begin{figure}[htbp]
\centerline{\epsfxsize=6.0in \epsfbox{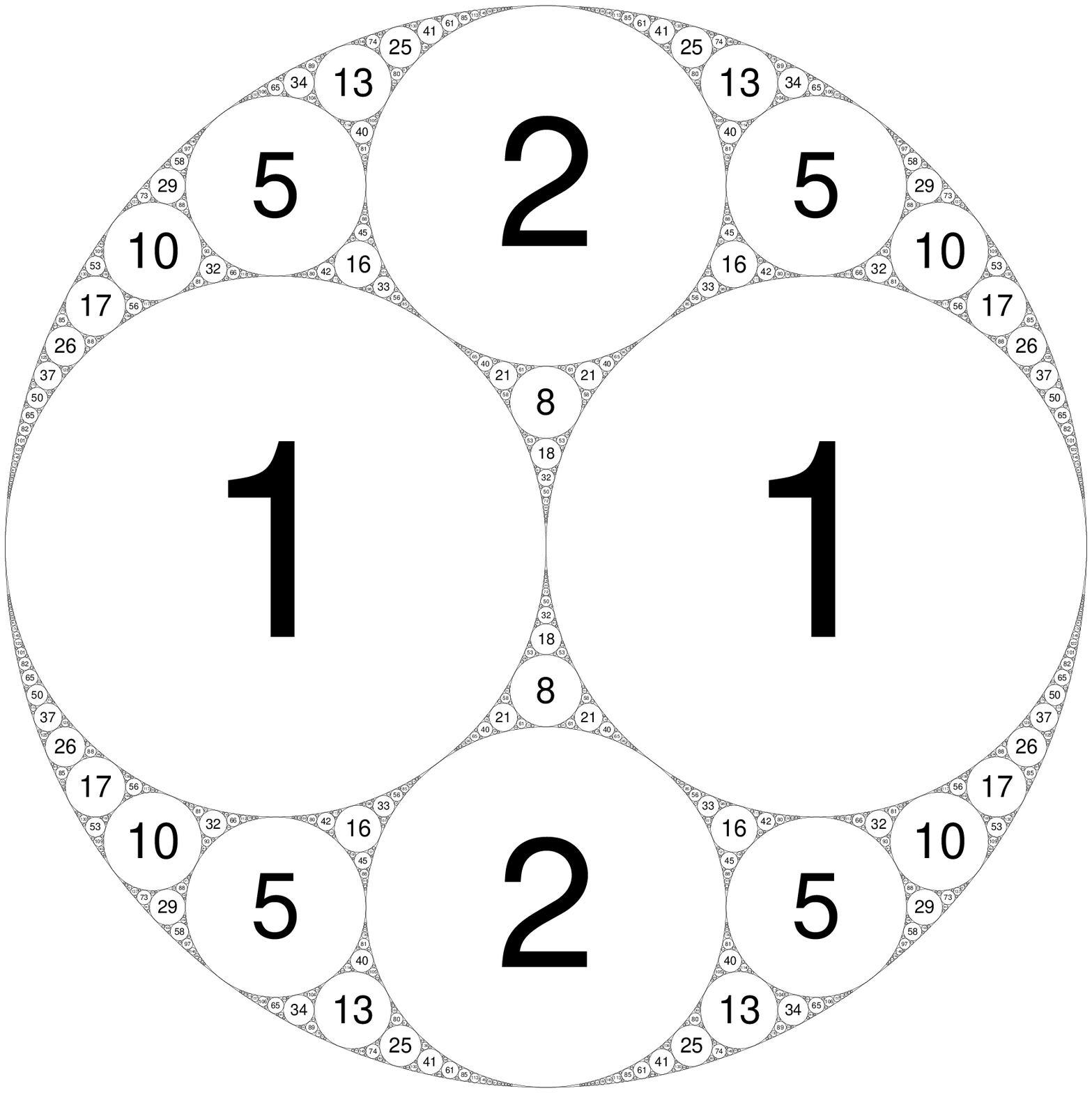}}
\caption{A spherical Apollonian packing}\label{fig6}
\end{figure}

Similarly, in the spherical packing associated to the standard Euclidean
packing in Figure 3, the $\cot \alpha$'s are all integers, different
from both the Euclidean and hyperbolic cases, starting from
$(0, 1, 1, 2).$ See Figure~\ref{fig6}.

One may notice interesting  numerical relations among the integers
in these three packings. Consider a ``loxodromic sequence''
of spheres as studied in Coxeter~\cite{Cox68b}, \cite{Cox97},
where each sphere is produced by reflection in the
largest sphere of the preceding Descartes configuration,
For the ``curvatures'' one  obtains for  the Euclidean packing 
the infinite sequence
$E: (-1, 2, 2, 3, 15, 38,...)$, for the spherical packing
$S: (0, 1, 1,  2, 8, 21, ...)$ and for the hyperbolic packing
$H: (-1, 1, 1, 1, 7, 17, ...)$. Note that $S + H = E,$ since
this is so for the initial values, and each sequence satisfies
the same fourth order linear recurrence relation, which is 
$x_{n + 1} = 2x_n + 2x_{n-1} + 2x_{n-2} - x_{n-3},$
by \eqn{201}.

In the Euclidean case, there are infinitely
many different kinds of Apollonian packings having integer curvatures
for all circles, see \cite{GLMWY2}. The same occurs for both 
hyperbolic and spherical 
Apollonian circle packings.
In the hyperbolic case  we 
also include among such integer hyperbolic circle
packings some packings which are ``virtual packings''.
Figure~\ref{fig5} is generated from the basic configuration having coth's
$(-1, 1, 1, 1)$,
and the next simplest hyperbolic case is $(-2, 3, 5, 6)$.

The Apollonian construction works also in higher dimensions,
but gives sphere packings only in dimensions two and three;  
in dimensions four and higher
we do not get proper packings; after several steps 
the spheres will overlap, see Boyd~\cite{Bo73}.
 However ``Apollonian sphere ensembles''
continue to exist 
in all dimensions as collections of Descartes configurations,
see \cite{GLMWY13}. 

There is a  considerable amount of 
mathematics devoted to circle packings;
Kenneth Stephenson's~\cite{Ste99} bibliography of circle-packing papers
lists over 90 papers since 1990. 
For further relations of
Apollonian packings and the relation of
integer Apollonian circle packings to the
integer Lorentz group $O(1, 3 , \zz)$,
see our series of papers with Ron Graham and Catherine Yan
\cite{GLMWY11},\cite{GLMWY12},\cite{GLMWY13},
 \cite{GLMWY2}  and S\"{o}derberg~\cite{So92}.

%
%

\section{Conclusion}

We have extended the Descartes circle theorem, well known
for  $n$-dimensional Euclidean space, to $n$-dimensional
spherical and hyperbolic space.
We presented matrix generalizations of the Descartes circle
theorem which
characterize Descartes configurations in all three geometries, 
and which required for their formulation the use of  
a particular coordinate system in each of these
geometries. Mauldon~\cite{Mau62} generalized the Soddy-Gossett
theorem in all three geometries to apply to sets of 
$n + 2$ equally inclined
spheres, as measured by an inclination parameter $\gamma$, 
with $\gamma = -1$ for touching spheres; our matrix theorems
can be extended to the case of arbitrary $\gamma$ as well.

Interestingly,  there are one-dimensional analogues of all these
theorems. For the
Euclidean case in 
one dimension a ``circle''
consists of two points bounding an interval, and two
``circles'' are tangent if they have one point in common.
The one-dimensional Descartes form is
\beql{701}
\bQ_1 :=  I_3 -  \bo_3 \bo_3^T= \left[ \begin{array}{ccc}
             0 & -1 & -1  \\
            -1 &  0 & -1   \\
            -1 & -1 &  0 \end{array} \right].
\eeq
A {\em one-dimensional Euclidean Descartes configuration }
consists of two touching intervals, 
and a third ``interval'' which is the complement
of their union, so that the three intervals cover the line $\rr.$
Call the third ``interval'' the {\em infinite interval,} 
and its  ``length'' is defined to be the negative of
the length of its complement,
which is the union of the first two intervals. The radius
is half the ``length.'' 
The radii $r_1, r_2, r_3,$ of the three intervals then satisfy
$$ r_1 + r_2 + r_3 = 0, $$
which is equivalent to the Descartes relation
\beql{702}
Q_1(\frac{1}{r_1}, \frac{1}{r_2}, \frac{1}{r_3}) = 
- \frac{2}{r_1r_2} - \frac{2}{r_1r_3} - \frac{2}{r_2r_3}  = 0.
\eeq
The value of ``curvature$\times$center'' of the infinite
interval is defined as being equal to the ``curvature$\times$center''
of the finite interval obtained  by reflection sending $x \to \frac{1}{x}.$
This describes  a positively oriented Descartes configuration; a negatively 
oriented one is obtained by reversing all signs. 
One can now define a $3 \times 3$ augmented matrix $\bW$ exactly as in
the augmented Euclidean Descartes theorem, and one finds that
\beql{703} 
\bW^T \bQ_1 \bW = \left[ \begin{array}{ccc}
             0 & -4 &  0  \\
            -4 &  0 &  0   \\
             0  &   0 &  2 \end{array} \right].
\eeq
Conversely, every solution $\bW$ to this equation corresponds to
a one-dimensional Descartes configuration.
The is even a notion of  Apollonian packing in dimension $n = 1$,
but it consists of
a single Descartes configuration! This holds
because there is only
a single circle tangent to a pair of tangent 
one-dimensional circles. 
That is,  the Descartes
equation \eqn{702} is linear in each curvature variable 
$a_i = \frac{1}{r_i}$ separately, 
instead of quadratic, hence the
reflection operation which generates  
new circles to add to the Apollonian packing in
dimensions $n \geq 2$ does not exist. Finally,
there are one-dimensional
spherical and hyperbolic analogues of these results,
defined via   \eqn{801}, taking $n = 1$. 
They can be established by stereographic
projection.

The main results in this paper are theorems in 
{\em inversive geometry}, as described in
Wilker~\cite{Wi81}, also Alexander~\cite{A67}
and Schwerdtfeger~\cite{Sch79}.
Inversive geometry is the
geometry that preserves spheres and their incidences,
which consists of the study of
geometric properties preserved by the group M\"{o}b$(n)$ of
conformal transformations of 
the space $\hat{\rr}^n = \rr^n \cup \{\infty\} \approx S^n.$
The set  of
Descartes configurations form a single orbit under the
action of the conformal group, and this group
appears~\footnote{The conformal group is isomorphic to a subgroup
of index $2$ in $Aut(Q_n)$, and we introduced
oriented Descartes configurations to keep track of the two cosets.}
 in our results as the (real) automorphism group
$$Aut(Q_n) := \{ \bN : ~ \bN^T \bQ_n \bN = \bQ_n \}$$
of the Descartes quadratic form $Q_n$, which is a Lie
group isomorphic to $O(n + 1, 1),$ (see
Wilker~\cite[Corollary p. 390 ]{Wi81} for the isomorphism),
and  the three generalized Descartes theorems given here 
 are invariant under the action of $Aut(Q_n).$
Our results prompt the question: Is 
there a ``natural'' characterization of the global coordinate
systems used in these geometries which yields the generalized
Descartes circle theorems in the matrix form presented here?

{\tt
\begin{tabular}{lll}
email: 
& jcl@research.att.com \\
  & clm@research.att.com \\
 & allan@research.att.com \\
\end{tabular}
 }

\end{document}